# LINDENMAYER SYSTEMS AND PRIMES

## About integer sequences containing unexpectedly many – or unexpectedly few – primes

by Andrei Vieru

**Abstract**
We study the surprising discrepancy between the number of primes corresponding, respectively, to A and to B in the infinite word ABBABAABBAABABBA... engendered by one of the simplest L-systems. We formulate a conjecture concerning the rate of growth of this discrepancy – which seems to tend to $e$ for every two sufficiently high consecutive even rank iterates of the L-system. We also submit a series of conjectures about a strange periodicity that appears when we compare the sequence of primes and the sequence of letters produced by this L-system (which, as everybody knows, are both aperiodic). Mysteriously, the periodicity is still obvious when we consider the sequences of integers – and, therefore, of primes – of the form $4x+1$ and of the form $4x+3$ ; strangely, it disappears when we consider only integers of the form $6x+1$ and of the form $6x+5$. Proofs or invalidations, comparisons with primes of other forms and/or with other Lindenmayer systems are proposed as unsolved or not yet studied problems to the readers.

## Introduction

Let $\Phi$ be the Lindenmayer system based on the alphabet {A, B} and the rules A→AB and B→BA, to be applied simultaneously to all letters of a given word. We shall call the words A, AB, ABBA, ABBABAAB, ABBABAABBAABABBA,... the $\Phi$ L-System's iterates of order 0, 1, 2, 3, 4,... respectively. In this L-system, any iterate of order $k$ will be designated by $<A,B>^k$. Replacing A by B and B by A in $<A,B>^k$ and pasting the result at the end of $<A,B>^k$, we'll obtain $<A,B>^{k+1}$. It is also well known that none of the iterates, including $<A,B>^\infty$ may be divided into periods. Every iterate may also be obtained from the next one by suppressing letters of even rank. More generally, suppressing in $<A,B>^k$ all letters of rank $2^m+n$ ($n \neq 0$, $n < 2^m$), i. e. keeping only letters of rank $2^m$, we will obtain $<A,B>^{k-m}$.
$<A,B>^\infty$ – which, as a matter of fact, contains an *uncountable* infinity of letters – is an invariant in regard to all operations described above.

Let us consider the iterate of infinite order $<A, B>^\infty$. Let S be one of its finite sub-sequences. Let $|S|$ be its total number of letters, while $|S_A|$ (respectively $|S_B|$) will designate the number of the occurrences of A in S (respectively of B in S).

Clearly, we'll always have $||S_A| - |S_B|| \leq 2$

If $|S|$ is odd, then $||S_A| - |S_B|| = 1$

If $|S|$ is even, either $|S_A| - |S_B| = 0$ or $||S_A| - |S_B|| = 2$

In particular, if $|S|$ is even and if S contains the beginning of some iterate, then $|S_A| - |S_B| = 0$

These considerations enable us to state that the densities of A and B in any iterate – in fact, in any beginning of any iterate – are equal. Moreover, the A and B densities are equal among letters of odd (and of even) rank.

## 1. Comparison with (prime) integers

Now let us fully write the beginning of some $<A,B>^k$ in order to clearly show the rank of each letter:

A  B  B  A  B  A  A  B  B  A  A  B  A  B  B  A  B  A  A  B  A …
↓  ↓  ↓  ↓  ↓  ↓  ↓  ↓  ↓  ↓   ↓   ↓   ↓   ↓   ↓   ↓   ↓   ↓   ↓   ↓   ↓
1  2  3  4  5  6  7  8  9  10 11 12 13 14 15 16 17 18 19 20 21 …

For every iterate of order *k*, one can ask: how many A and how many B are of prime rank? Intuition might have already whispered that the probability for a prime integer to correspond to an A equals the probability to correspond to a B. In other words, intuition says we will find a series of A and B of prime rank very much resembling to the result of a 'repeated experiment with two equiprobable outcomes'. If this intuition were correct, it would imply that, replacing A by –1 and B by +1, then calculating the sum of the outcomes after every step[1], we'll have 0 ('the standard predictable value') as many times as we wish, and that the set of these 'partial' sums is not bounded: it can be as big in absolute value as we wish (so-called 'probabilistic law of the pendulum'). In respect to 'the standard predictable value', the 'predictable error' is measured by means of the 'dispersion unit', which is, in our case, $n^{½}/2$, where *n* represents the number of steps after which we check the sum (or 'error').

In reality, although in the VERY long run the probability to find an A (or a B) when looking to a letter of prime rank seems to approach 0.5, the series of results of this 'repeated experiment with two (supposedly) equiprobable outcomes' does not resemble at all to the classical case of a 'repeated throw of a coin'.

---
[1] after every prime

Let's designate by $|Pr(S_A)|$ (respectively by $|Pr(S_B)|$) the number of A (respectively of B), of prime rank in S. In particular, when $S = <A, B>^k$ we'll write, when needed $|Pr(<A,B>^k_A)|$ and $|Pr(<A,B>^k_B)|$, or, shortly, $|Pr(A^k)|$ and $|Pr(B^k)|$

Let's designate by $|Pr(A(n))|$, respectively by $|Pr(B(n))|$, the number of A, respectively B, of prime rank, encountered among the first $n$ letters of $<A,B>^\infty$.

**CONJECTURE 1**
$|Pr(A^k)| - |Pr(B^k)| \to \infty$ as $k \to \infty$
$|Pr(A^k)|/|Pr(B^k)| \to 1$ as $k \to \infty$

In fact, $\min\{x \in Z \mid \exists n \in N \; x=|Pr(A(n))| - |Pr(B(n))|\} = -3$
This minimum is reached for $n = 5$ and never again[2].

**CONJECTURE 2**
$\forall m \in Z \quad m = |Pr(A(x))| - |Pr(B(x))|$ has only a finite number of solutions, or no solutions at all. (In the case of the 'repeatedly thrown coin' we normally have infinitely many 'solutions' for every positive or negative integer.)

## 2. Rate of growth of the discrepancy
Here are the $m$ values – i. e. $|Pr(A^k)| - |Pr(B^k)|$ – for $k = 10$, $k = 11$, $k = 12$, $k = 13$, $k = 14$, $k = 15$, $k = 16$, $k = 17$, $k = 18$, $k = 19$, $k = 20$, $k = 21$, $k = 22$, $k = 23$, $k = 24$ :

*38, 53, 116, 141, 274, 333, 824, 829, 1749, 1897, 5061, 5220, 13515, 14969, 36721…*

Let's write, shorter, $D_1(n) = |Pr(A(n))| - |Pr(B(n))|$
As one can see, generally we have $2 < D_1(2^{k+2})/ D_1(2^k) < 3$
In particular $36721/13515 = 2.717055…$ which is already surprisingly near $e$… Although there is not yet enough evidence, I take the risk to propose the following

**CONJECTURE 3**
For a chosen $k \geq 3$, $\lim_{j \to \infty} \ln D_1(2^{2k + 2j}) - \ln D_1(2^{2k}) - j = 0$

---

[2] In the case of 'normal' series of experiments with two equiprobable outcomes, this minimum do not exist.

or $\quad \lim_{\substack{j \to \infty \\ k \to \infty}} [D_1(2^{2k+2j}) / D_1(2^{2k})]^{1/j} = e$

## 3. Shifts and periodicity

When we reckon the letters $<A,B>^\infty$ is made of, we can try to see what happens if, instead of starting with 1, we begin counting with some arbitrarily chosen (positive or negative) integer:

A   B   B   A   B   A   A   B   etc.
↓   ↓   ↓   ↓   ↓   ↓   ↓   ↓
$z$  $z+1$  $z+2$  $z+3$  $z+4$  $z+5$  $z+6$  $z+7$   etc.

We'll call this way of counting the ranks of each letter 'count with shift $z$'. Considering the first $n$ letters in some $<A,B>^k$ iterate and counting them with shift $z$, let's designate by $|Pr(A_z(n))|$, respectively by $|Pr(B_z(n))|$, the number of A – respectively of B – that correspond to a prime. (Thus, we should write from now on $|Pr(B_1(n))|$ instead of $|Pr(B(n))|$.)

**CONJECTURE 4**

For every $z$ 'not too big in respect to $k$':

$|Pr_{6z}(A^k)|-|Pr_{6z}(B^k)| \to -\infty$ and $|Pr_{6z}(A^k)|/|Pr_{6z}(B^k)| \to 1$ as $k \to \infty$ \qquad (0)

$|Pr_{6z+1}(A^k)|-|Pr_{6z+1}(B^k)| \to +\infty$ and $|Pr_{6z+1}(A^k)|/|Pr_{6z+1}(B^k)| \to 1$ as $k \to \infty$ \qquad (1)

$|Pr_{6z+2}(A^k)|-|Pr_{6z+2}(B^k)| \to +\infty$ and $|Pr_{6z+2}(A^k)|/|Pr_{6z+2}(B^k)| \to 1$ as $k \to \infty$ \qquad (2)

$|Pr_{6z+3}(A^k)|-|Pr_{6z+3}(B^k)| \to -\infty$ and $|Pr_{6z+3}(A^k)|/|Pr_{6z+3}(B^k)| \to 1$ as $k \to \infty$ \qquad (3)

$|Pr_{6z+4}(A^k)|-|Pr_{6z+4}(B^k)| \to -\infty$ and $|Pr_{6z+4}(A^k)|/|Pr_{6z+4}(B^k)| \to 1$ as $k \to \infty$ \qquad (4)

$|Pr_{6z+5}(A^k)|-|Pr_{6z+5}(B^k)| \to +\infty$ and $|Pr_{6z+5}(A^k)|/|Pr_{6z+5}(B^k)| \to 1$ as $k \to \infty$ \qquad (5)

**CONJECTURE 5**

In (4) and (5), convergence to $-\infty$ and $+\infty$ are much slower than in (0), (1), (2) and (3), while convergence to 1 is faster.
For every $k$ and for every $z$ 'not too big in respect to $k$', $|Pr_{6z+4}(A^{2k-1})| - |Pr_{6z+4}(B^{2k-1})| = |Pr_{6z+4}(A^{2k})| - |Pr_{6z+4}(B^{2k})| = |Pr_{6z+5}(B^{2k-1})| - |Pr_{6z+5}(A^{2k-1})| = |Pr_{6z+5}(B^{2k})| - |Pr_{6z+5}(A^{2k})|$

Commentary to conjectures 4 and 5:

Scanning, for all $n$, $D_4(n) = |Pr(A_4(n))| - |Pr(B_4(n))|$ and $D_5(n) = |Pr(A_5(n))| - |Pr(B_5(n))|$ it is not yet sure neither that $D_4(n)$ ($n \neq 2^k$) is positively bounded (has a maximum) nor that $D_5(n)$ ($n \neq 2^k$) is negatively bounded (has a minimum).

**CONJECTURE 6**
Now, one can ask what happens if we set up a relationship between iterates of $\Phi$ and integers of the form $4n + 1$ or $4n + 3$:
A   B   B   A   B   A   A   B …
1   5   9   13  17  21  25  29 …
3   7   11  15  19  23  27  31 …

Well, we'll find the same discrepancy when counting primes of the form $4n + 1$ (or, in a separate experiment, of the form $4n + 3$) corresponding to A and to B. We'll find the same periodicity of 6, when we'll count primes while introducing 'shifts', and we can formulate a set of conjectures analogous to (0), (1), (2), (3), (4) and (5) in 'Conjecture 3'.

**CONJECTURE 7**
One can ask what happens if we establish a relationship between $\Phi$'s iterates and integers of form $6n + 1$ or $6n + 5$:
A   B   B   A   B   A   A   B   B…
1   7   13  19  25  31  37  43  49…
5   11  17  23  29  35  41  47  53…

Whatever shift in counting we might set up, for every iterate $<A,B>^k$, the discrepancy between the number of primes of the form $6n + 1$ (or, in a distinct experiment, of the form $6n + 5$) corresponding to A and to B becomes so little that it is not anymore possible to formulate any conjecture that would have establish any significant difference between this experiment and a 'repeated throw of a coin'. Needless to say, there is not anymore any periodicity of 6 (nor any other periodicity) similar to the periodicity 'Conjecture 4' is based on.

**Conclusion**
So far, we haven't yet studied $\Phi$'s iterates in relationship with integers of the form $pn+a$ (where $p$ is prime). Nor have we studied L-systems based, for example, on an alphabet $\{A_1, A_2,…, A_n\}$ and on rules of the form $A_i \rightarrow A_{i(1)}A_{i(2)}…A_{i(n)}$ such as there is a group where $\forall i\ \forall j\ A_i + A_j = A_{i(j)}$ (Example: alphabet $\{A, B, C\}$ and rules A→ABC  B→BCA  and C→CAB.) These comparisons might be interesting, although the results would be even more difficult to interpret.

      Andrei Vieru

**Acknowledgments**
We express our gratitude to Vlad Vieru for skilful computer programming and to Robert Vinograd for useful discussions.